\documentclass[journal]{IEEEtran}

\ifCLASSINFOpdf
\else
   \usepackage[dvips]{graphicx}
\fi
\usepackage{url}

\usepackage{graphicx}
\usepackage{amsmath,amsfonts,amssymb,amsbsy, amsthm,ae,aecompl}
\usepackage{algorithm,algpseudocode}
\usepackage[utf8]{inputenc}
\usepackage[english]{babel}
\usepackage{bm}
\usepackage{color,soul}
\usepackage[noadjust]{cite}
\usepackage{epsfig}
\usepackage{enumerate}
\usepackage{float}
\usepackage{fancyhdr}
\usepackage[T1]{fontenc}
\usepackage[acronym,toc,shortcuts]{glossaries}
 \usepackage[caption=false,font=footnotesize]{subfig}
 \usepackage{soul}
\usepackage{lipsum}
\usepackage{dsfont}
\usepackage{transparent}
\usepackage{tikz,pgfplots}
\usepackage{times}
\usepackage{verbatim}
\usepackage[all]{xy}
\usepackage[colorinlistoftodos]{todonotes}
\usepackage{xcolor}

\makeglossaries
\newacronym[sort=ell]{test}{LED}{limiting eigenvalue distribution}
\newacronym{RGG}{RGG}{random geometric graph}
\newacronym{DGG}{DGG}{deterministic geometric graph with nodes in a grid}
\newacronym{RGGs}{RGGs}{random geometric graphs}
\newacronym{ESDF}{ESDF}{ empirical spectral distribution function}
\newacronym{ESDFs}{ESDFs}{ empirical spectral distribution functions}
\newacronym{ER}{ER}{Erd\"{o}s-Rényi}
\newacronym{LSD}{LSD}{limiting spectral distribution}

\newtheorem{theorem}{Theorem}

\newtheorem{lemma}{Lemma}

\begin{document}
\title{Spectral bounds of the regularized normalized Laplacian for random geometric graphs}
\author{\IEEEauthorblockN{Mounia Hamidouche$^{\star}$, Laura Cottatellucci$^{\dagger}$, Konstantin Avrachenkov$^{ \diamond}$}
\IEEEauthorblockA{
				\vspace{0.2cm} \\			
		{		\small  $^{\star}$ Departement of communication systems, EURECOM, France \\	
		 $^{\dagger}$ Department of Electrical, Electronics, and Communication Engineering, FAU, Germany \\	
				$^{\diamond}$ INRIA Sophia Antipolis, France.	 }
		} 
\thanks{This research was funded by the French Government through the Investments for the Future Program with Reference: Labex UCN@Sophia-UDCBWN.
 }}
\maketitle
\begin{abstract}
In this work, we study the spectrum of the regularized normalized Laplacian for random geometric graphs (RGGs) in both the connectivity and thermodynamic regimes. We prove that the limiting eigenvalue distribution (LED) of the normalized Laplacian matrix for an RGG converges to the Dirac measure in one in the full range of the connectivity regime. In the thermodynamic regime, we propose an  approximation for the LED and we provide a bound on the Levy distance between the approximation and the actual distribution. In particular, we show that the LED of the  regularized normalized Laplacian matrix for an RGG can be approximated by the LED of the regularized normalized Laplacian for a deterministic geometric graph with nodes in a grid (DGG). Thereby, we obtain an explicit approximation of the eigenvalues in the thermodynamic regime.
\end{abstract}
\IEEEpeerreviewmaketitle
\section{Introduction}

\IEEEPARstart{T}{he} spectrum of the normalized Laplacian associated to \glspl{RGG} provides very useful insights on the structure and the dynamics of large complex networks in which  the geographical distance is a critical factor.

An \gls{RGG} is a graph with a finite set $\mathcal{X}_{n}$ of $n$ points $x_{1},...,x_{n},$ distributed uniformly and independently on a $d$-dimensional torus $\mathbb{T}^d \equiv [0, 1]^d$.  Given a geographical distance $r_{n} >0 $, we form a graph by connecting two points $x_{i}, x_{j} \in \mathcal{X}_{n}$ if the distance between them is at most $r_{n}$, i.e., $\mathrm{d}\left(x_{i}, x_{j} \right) \leq r_{n}$. 
 The maximum distance $r_{n}$ is a function of $n$ chosen such that $r_{n}\rightarrow 0$ when $n \rightarrow \infty$.  The average vertex degree in $G(\mathcal{X}_{n}, r_{n})$ is $a_{n}=\theta^{(d)} nr_{n}^d$, where $\theta^{(d)}$ denotes the volume of a $d$-dimensional unit hypersphere in $\mathbb{T}^d$ \cite{penrose2003random}. An \gls{RGG} presents different properties depending on the average vertex degree $a_{n}$. Two different scaling regimes for $a_{n}$ are of interest. The first one is the connectivity regime, in which $a_{n}$ scales as $\Omega(\log(n))$\footnote{The notation $f(n) =\Omega(g(n))$ indicates that $f(n)$ is bounded below by $g(n)$ asymptotically, i.e., $\exists K>0$ and $ n_{o} \in \mathbb{N}$ such that $\forall n > n_{0}$ $f(n) \geq K g(n)$.} \cite{penrose2003random}. The second scaling regime of interest is the thermodynamic regime in which $a_{n}\equiv \gamma$, for $\gamma >0$. 
 
 In \cite{rai2007spectrum}, the author shows that the \glspl{test} of the transition probability matrix of random walks in an \gls{RGG} and in a \gls{DGG} converge to the same limit in the connectivity regime. However, for $\epsilon >0$ the result in  \cite{rai2007spectrum} holds only when $a_{n}$ scales as $\Omega\left(\sqrt{n}\log^{\epsilon}(n) \right)$ for $d=1$, as $\Omega\left(\log^{\frac{3}{2}+\epsilon}(n)\right)$ for $d=2$, and as $\Omega\left(\log^{1+\epsilon}(n)\right)$ for $d\geq 3$. Motivated by the result in \cite{rai2007spectrum}, here we study the \gls{test} of the regularized normalized Laplacian for an \gls{RGG} in the full range of the connectivity regime, i.e, $a_{n}=\Omega(\log(n))$, and in the thermodynamic regime. To the best of our knowledge, an explicit expression for the \gls{test} of the normalized Laplacian for \glspl{RGG} is still not known in the full range of the scaling laws for the average vertex degree $a_{n}$ in the connectivity regime, nor in the thermodynamic regime.
 
Using the above defined \gls{RGG}, we show that the \glspl{test} of the normalized Laplacian for an \gls{RGG} and for a \gls{DGG} converge to same limit for any dimension $d \geq 1$ as $n \rightarrow \infty$ in the full range of the connectivity regime. Then, we show that their corresponding \glspl{test} converge to the Dirac measure in one as $n \rightarrow \infty$ and $d=1$. In the thermodynamic regime, we provide an analytical approximation for the eigenvalues of the regularized normalized Laplacian for the \gls{RGG} when $d=1$. Then, we numerically validate our theoretical results by comparing the simulated and the analytical spectum.

 \section{ Spectral Analysis}
 \label{Sec1}
 In order to study the \gls{test} of \glspl{RGG}, we consider a finite set $\mathcal{D}_{n}$ of $n$ grid points that are at the intersections of all parallel hyperplanes with separation $n^{-1/d}$. Then, we define a deterministic geometric graph in a grid $G(\mathcal{D}_{n}, r_{n})$ by connecting two points in $\mathcal{D}_{n}$ when the distance between them is at most $ r_{n}$. Given two nodes in $G(\mathcal{X}_{n}, r_{n})$ or in $G(\mathcal{D}_{n}, r_{n})$, we assume that there is at most one edges between them and there is no edge from a vertex to itself. We denote the degree of a vertex in $G(\mathcal{D}_{n}, r_{n})$ by $a'_{n}$, and in particular $a'_{n}\equiv\gamma'$ in the thermodynamic regime. Let $\mathbf{N}(x_{i})$ be the number of neighbors of a vertex $x_{i}$ in $G(\mathcal{X}_{n}, r_{n})$. Define $\hat{\mathcal{L}}(\mathcal{X}_{n})$ as the regularized normalized Laplacian matrix for $G(\mathcal{X}_{n}, r_{n})$ with entries,
\begin{equation} \hat{\mathcal{L}}(\mathcal{X}_{n})_{i j}= \delta_{ij}-\dfrac{\chi [x_{i} \sim x_{j}]+\frac{\alpha}{n}}{ \sqrt{(\mathbf{N}(x_{i})+\alpha)(\mathbf{N}(x_{j})+\alpha)}}.
 \end{equation}
The term $\chi[x_{i}\thicksim x_{j}]$ takes unit value when there is an edge between nodes $i$ and $j$ in $G(\mathcal{X}_{n}, r_{n})$ and zero otherwise. $\delta_{ij}$ is the Kronecker delta function and $\frac{\alpha}{n} \geq 0$. A similar definition holds for $\hat{\mathcal{L}}(\mathcal{D}_{n})_{i j}$, the regularized normalized Laplacian defined over the nodes in $G(\mathcal{D}_{n}, r_{n})$. The regularized version of the normalized Laplacian matrix is used to overcome the problem of singularities due to isolated vertices in the thermodynamic regime \cite{avrachenkov2010improving}. 

In the following Theorem {\ref{theorem}} and {\ref{corollary1:general2}}, we show that the empirical spectral distribution function $F^{ \hat{\mathcal{L}}({\mathcal{D}_{n}})}$ is a good approximation of the spectral distribution $F^{  \hat{\mathcal{L}}({\mathcal{X}_{n}})}$ when $n$ is large in both the connectivity and thermodynamic regimes using the Levy distance $L^{3}$ \cite{taylor2012introduction}.
\begin{theorem}
\label{theorem} 
In the connectivity regime, i.e., as $a_{n}=\Omega (\log(n))$, when $d \geq 1$, $a_{n} \geq 2d$, $\alpha \rightarrow 0$ and $n \rightarrow \infty$ then, for any $t>0$,
\begin{equation*}
\lim_{n \to\infty} P \left\lbrace L^{3} \left( F^{\hat{\mathcal{L}}(\mathcal{X}_{n})}, F^{\hat{\mathcal{L}}(\mathcal{D}_{n})}\right) > t \right\rbrace  = 0.
\end{equation*}
\end{theorem}
In particular, when $a_{n}=c\log{(n)}$, $c>24$ then, the \gls{test} of the normalized Laplacian for an \gls{RGG} converges to the \gls{test} of the normalized  Laplacian for a \gls{DGG} with rate of convergence $\mathcal{O}\left( 1/ n^{c/12-1}\right)$, and when $c\leq 24,$  the rate of convergence is $\mathcal{O}\left( 1/n\right)$.  For $\epsilon >0$ and  $a_{n} \geq \log^{1+\epsilon}{(n)}$,  the rate of convergence is $\mathcal{O}\left( 1/n^{(a_{n}/12 \log{(n)})-1} \right)$. An exponential rate of convergence $\mathcal{O}\left(n e^{-n/12}\right)$ holds when the graph is dense, i.e., $a_{n}$ scales as $\Omega(n)$.  Hence, the result given in Theorem \ref{theorem} generalizes the work in \cite{rai2007spectrum} and shows that the \glspl{test} of the normalized Laplacian of $G(\mathcal{X}_{n}, r_{n})$ and $G(\mathcal{D}_{n}, r_{n})$ converge to the same limit as $n \rightarrow \infty$ in the full range of the connectivity regime and for any chosen distance function.
\begin{theorem}
\label{corollary1:general2} 
In the thermodynamic regime, i.e., as $a_{n} \equiv \gamma$,  for $d \geq1$, $\gamma \geq 2d$ and for every $t > \dfrac{8 \gamma}{(\gamma'+\alpha)^2}$
\begin{equation*}
\lim_{n\to\infty} P \left\lbrace L^{3} \left( F^{  \hat{\mathcal{L}}(\mathcal{X}_{n})}, F^{  \hat{\mathcal{L}}(\mathcal{D}_{n})} \right) > t \right\rbrace =0.
\end{equation*}
\end{theorem}
From Theorem {\ref{corollary1:general2}} it is apparent that, as $n \rightarrow \infty$ and $\alpha \rightarrow 0$,  $F^{  \hat{\mathcal{L}}(\mathcal{D}_{n})}$ approximates $F^{  \hat{\mathcal{L}}(\mathcal{X}_{n})}$ with an error bound of $8 \gamma/\gamma^{'2}$.
Now, we provide the \gls{test} of the normalized Laplacian in the connectivity regime as $n \rightarrow \infty$. By utilizing Theorem {\ref{theorem}} and the expression for the eigenvalues of the adjacency matrix of a {\gls{DGG}} in {\cite{nyberg2014laplacian}}, we can state the following Lemma {\ref{dirac}}.
 \begin{lemma}
\label{dirac}
In the connectivity regime, for $d=1$, $\alpha \rightarrow 0$ and using the Euclidean distance, the eigenvalues $\lambda_{m}$ of $ \hat{\mathcal{L}}(\mathcal{D}_{n})$ are obtained from
\begin{equation}
\lambda_{m}=1 -\dfrac{1}{a'_{n}}   \left(  \dfrac{\sin(\frac{m \pi}{n}(a'_{n}+1))}{\sin(\frac{m \pi}{n})}  - 1 \right),
\label{eig4}
\end{equation}
where $m$ $\in$ $\lbrace 0,...n-1 \rbrace$, $a'_{n}=2\left\lfloor n r_{n} \right\rfloor$ and $ \left\lfloor x \right\rfloor$ is the integer part, i.e.,  the greatest integer less than or equal to $x$. Then, as  $n\rightarrow \infty,$ the {\gls{test}} of $\hat{\mathcal{L}}(\mathcal{D}_{n})$ converges to the Dirac measure in one.
\end{lemma}
Now, we give an  approximation of the eigenvalues for the regularized normalized Laplacian of an \gls{RGG} in the thermodynamic regime as $n \rightarrow \infty$.

\begin{lemma}
\label{eig2}
 In the thermodynamic regime, for $d = 1$ and using the Euclidean distance, the asymptotic eigenvalues of $ \hat{\mathcal{L}}(\mathcal{D}_{n})$, as $n \rightarrow \infty$ are given by 
 \begin{equation} \label{eig}
\begin{split}
\lambda_{w}  & = 1- \dfrac{1}{(\gamma'+\alpha)}  \dfrac{\sin(w(\gamma'+1))}{\sin(w)}  +\dfrac{1-\alpha\delta_{w}}{(\gamma'+\alpha)},
\end{split}
\end{equation}
where $w \in [0, \pi]$, $\gamma'=2\left\lfloor \gamma \right\rfloor$ for $\gamma \geq 2$ and $\delta_{w}=1$ when $w=0$ otherwise $\delta_{w}=0.$ 
\end{lemma}
In particular, when $\alpha \rightarrow 0$ then, (\ref{eig}) reduces to 
\begin{equation}
\lambda_{w}=1 -\dfrac{1}{\gamma'}   \left(  \dfrac{\sin(w(\gamma'+1))}{\sin(w)}  - 1 \right).
\end{equation}

\section{Numerical Evaluation}

\begin{figure}[h!]
\subfloat[Connectivity regime for different $n$ and $r_{n}=\log^{\frac{3}{2}}(n)/n$.]{
\includegraphics[height=5.3cm, width=\columnwidth]{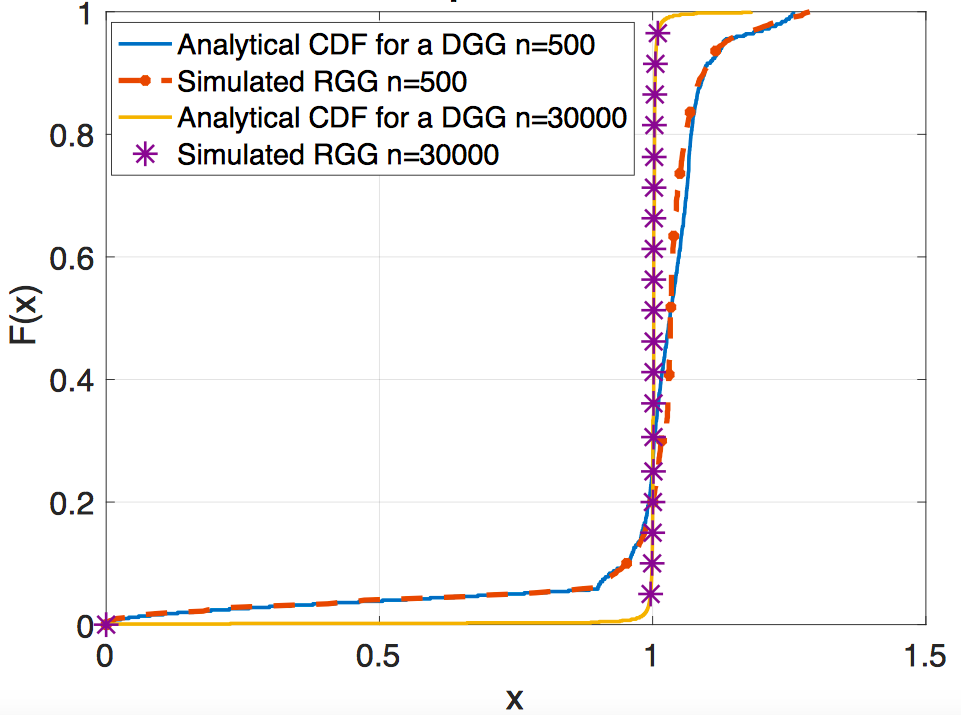}}
\hfill
\subfloat[Thermodynamic regime for $\gamma=12$ and $r_{n}=12/n$.]{\includegraphics[height=5.3cm, width=\columnwidth]{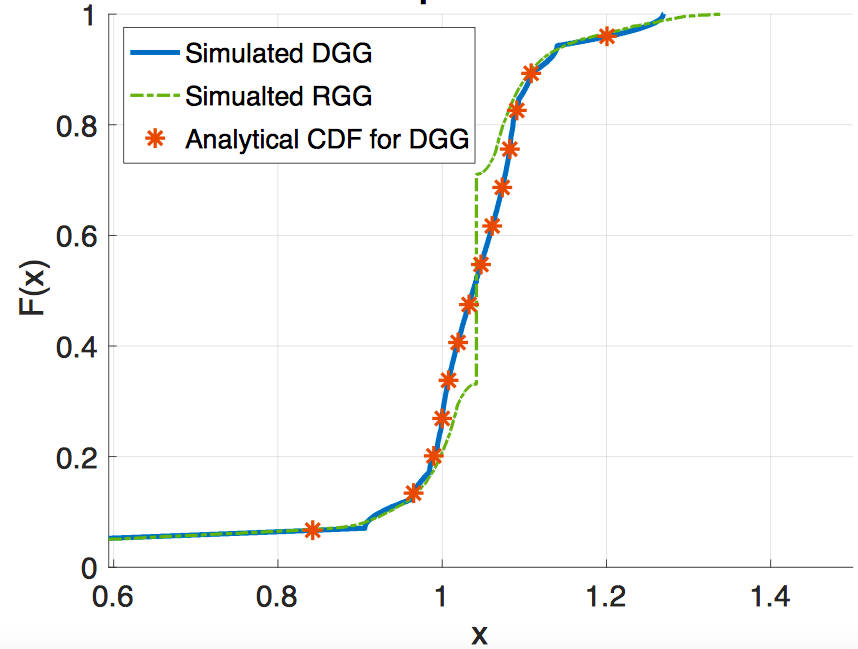}}
\caption{ Spectral distribution of an RGG and a DGG.}
\label{fig2}
\end{figure}

We validate our analytical results obtained in Section \ref{Sec1} on the spectrum of the regularized normalized Laplacian of \glspl{RGG} by numerical computations. In Fig. \ref{fig2}.a we compare the spectral distribution of a {\gls{DGG}} (continuous line) with the one of an {\gls{RGG}} with increasing the number of nodes $n$ (dashed line for $n=500$ and star markers for $n=30000$) in the connectivity regime. We notice that the curves corresponding to the \gls{RGG} and the \gls{DGG} match very well when $n$ is large. This confirms the result given in Theorem \ref{theorem}. It appears that by increasing $n$, the eigenvalue distribution converges to the Dirac measure in one as shown in Lemma \ref{dirac}. Fig. \ref{fig2}.b  illustrates the empirical spectral distribution of a simulated regularized normalized Laplacian in the thermodynamic regime for a realization of an {\gls{RGG}} and {\gls{DGG}} with $n=30000$ nodes and $\alpha =0.001$. The analytical distribution is obtained using Lemma {\ref{eig2}}. The gap that appears between the eigenvalue distributions of the \gls{RGG} and the \gls{DGG} is quantified by the error bound in Theorem \ref{corollary1:general2}.

\bibliographystyle{IEEEtran}
\bibliography{references}

\end{document}